%21-9-99

\magnification 1200

\input amstex
\input amsppt.sty
\vsize 18.7cm
\hsize 13.5cm

%\input nume2.tex

%%%%%%%%%%%%%%%%% EQUAZIONI CON NOMI SIMBOLICI
%%%
%%% Per assegnare un nome simbolico ad una equazione basta
%%% scrivere \Eq(...) o, in \eqalignno, \eq(...) o,
%%% nelle appendici, \Eqa(...) o \eqa(...);
%%% dentro le parentesi e al posto di ... si puo' scrivere qualsiasi commento;
%%% per avere i nomi simbolici segnati a sinistra delle formule si deve
%%% dichiarare il documento come bozza, iniziando il testo con
%%% \BOZZA. Sinonimi: \Eq,\EQ,\EQS; \eq,\eqs; \Eqa,\Eqas;\eqa,\eqas.
%%% All' inizio di ogni paragrafo si devono definire il
%%% numero del paragrafo e della prima formula dichiarando
%%% \numsec=... \numfor=...  (brevetto Eckmannn).
%%% Si possono citare formule seguenti; le corrispondenze fra nomi
%%% simbolici e numeri effettivi sono memorizzate nel file \jobname.aux, che
%%% viene letto all'inizio, se gia' presente. E' possibile citare anche
%%% formule che appaiono in altri file, purche' sia presente il
%%% corrispondente file .aux; basta includere all'inizio l'istruzione
%%%           \include{nomefile}
%%%
%%%%%%%%%%%%%%%%%%%%%%%%%%%%%%%%%%%%%%%%%%%%%%%%%%%%%%%%%%%%%%%

\global\newcount\numsec
\global\newcount\numfor
\global\newcount\numtheo
\global\advance\numtheo by 1

\def\senondefinito#1{\expandafter\ifx\csname#1\endcsname\relax}

\def\SIA #1,#2,#3 {\senondefinito{#1#2}%
\expandafter\xdef\csname #1#2\endcsname{#3}\else
\write16{???? ma #1,#2 e' gia' stato definito !!!!} \fi}

\def\etichetta(#1){(\veroparagrafo.\veraformula)%
\SIA e,#1,(\veroparagrafo.\veraformula) %
\global\advance\numfor by 1%
\write15{\string\FU (#1){\equ(#1)}}%
\write16{ EQ #1 ==> \equ(#1) }}

\def\letichetta(#1){\veroparagrafo.\verotheo
\SIA e,#1,{\veroparagrafo.\verotheo}
\global\advance\numtheo by 1
 \write15{\string\FU (#1){\equ(#1)}}
 \write16{ Sta \equ(#1) == #1 }}

\def\tetichetta(#1){\veroparagrafo.\veraformula %%%%copy four lines
\SIA e,#1,{(\veroparagrafo.\veraformula)}
\global\advance\numfor by 1
\write15{\string\FU (#1){\equ(#1)}}
\write16{ tag #1 ==> \equ(#1)}}

\def\FU(#1)#2{\SIA fu,#1,#2 }

\def\etichettaa(#1){(A\veroparagrafo.\veraformula)%
\SIA e,#1,(A\veroparagrafo.\veraformula) %
\global\advance\numfor by 1%
\write15{\string\FU (#1){\equ(#1)}}%
\write16{ EQ #1 ==> \equ(#1) }}

\def\BOZZA{
\def\alato(##1){%
 {\rlap{\kern-\hsize\kern-1.4truecm{$\scriptstyle##1$}}}}%
\def\aolado(##1){%
 {%\vtop to \profonditastruttura
{%\baselineskip
 %\profonditastruttura\vss
 \rlap{\kern-1.4truecm{$\scriptstyle##1$}}}}}
}

\def\alato(#1){}
\def\aolado(#1){}

\def\veroparagrafo{\number\numsec}
\def\veraformula{\number\numfor}
\def\verotheo{\number\numtheo}

\def\Eq(#1){\eqno{\etichetta(#1)\alato(#1)}}
\def\eq(#1){\etichetta(#1)\alato(#1)}
\def\teq(#1){\tag{\aolado(#1)\tetichetta(#1)\alato(#1)}}%%%%%this line for \tag
\def\Eqa(#1){\eqno{\etichettaa(#1)\alato(#1)}}
\def\eqa(#1){\etichettaa(#1)\alato(#1)}
\def\eqv(#1){\senondefinito{fu#1}$\clubsuit$#1
\write16{#1 non e' (ancora) definito}%
\else\csname fu#1\endcsname\fi}
\def\equ(#1){\senondefinito{e#1}\eqv(#1)\else\csname e#1\endcsname\fi}

%next six lines by paf (no responsibilities taken)
\def\Lemma(#1){\aolado(#1)Lemma \letichetta(#1)}%
\def\Theorem(#1){{\aolado(#1)Theorem \letichetta(#1)}}%
\def\Proposition(#1){\aolado(#1){Proposition \letichetta(#1)}}%
\def\Corollary(#1){{\aolado(#1)Corollary \letichetta(#1)}}%
\def\Remark(#1){{\noindent\aolado(#1){\bf Remark \letichetta(#1)}}}%
\def\Definition(#1){{\noindent\aolado(#1){\bf Definition \letichetta(#1)$\!\!$\hskip-1.6truemm}}}
\def\Example(#1){\aolado(#1) Example \letichetta(#1)$\!\!$\hskip-1.6truemm}

\def\include#1{
\openin13=#1.aux \ifeof13 \relax \else
\input #1.aux \closein13 \fi}

\openin14=\jobname.aux \ifeof14 \relax \else
\input \jobname.aux \closein14 \fi
\openout15=\jobname.aux

%%% Local Variables: 
%%% mode: plain-tex
%%% TeX-master: t
%%% End: 

\parskip 2mm

\def\E{{\Bbb E}}
\def\P{{\Bbb P}}
\def\R{{\Bbb R}}
\def\Z{{\Bbb Z}}

\def\sqr{\vcenter{
         \hrule height.1mm
         \hbox{\vrule width.1mm height2.2mm\kern2.18mm\vrule width.1mm}
         \hrule height.1mm}}                  % This is a slimmer sqr.
\def\square{\ifmmode\sqr\else{$\sqr$}\fi}
\def\one{{\bold 1}\hskip-.5mm}

%\def\MM{{\bf m}}

%\BOZZA

\topmatter

\title
Rate of convergence to equilibrium of symmetric simple exclusion processes
\endtitle

\rightheadtext{Convergence to equilibrium of the exclusion process}

\author  P. A. Ferrari, A. Galves, C. Landim \endauthor

\affil{
Universidade de S\~{a}o Paulo \\
Universidade de S\~{a}o Paulo \\
IMPA and CNRS UPRES-A 6085, Rouen 
}  \endaffil

\address{IME USP,
Caixa Postal 66281,
 05315-970 - S\~{a}o Paulo,
 BRAZIL
} \endaddress
\email pablo\@ime.usp.br  \endemail

\address{IME USP,
Caixa Postal 66281,
 05315-970 - S\~{a}o Paulo,
 BRAZIL
} \endaddress
\email galves\@ime.usp.br  \endemail

\address{IMPA, Estrada Dona Castorina 110,
CEP 22460 Rio de Janeiro, Brasil and CNRS UPRES-A 6085,
Universit\'e de Rouen, 76128 Mont Saint Aignan, France.} \endaddress
\email  landim\@impa.br \endemail

\abstract We give bounds on the rate of convergence to equilibrium of
the symmetric simple exclusion process in $\Z^d$. Our results include
the existent results in the literature. We get better bounds and
larger class of initial states via a unified approach. The method
includes a comparison of the evolution of $n$ interacting particles
with $n$ independent ones along the whole time trajectory.
\endabstract

\keywords
interacting particle system, 
\endkeywords

\subjclass Primary 60K35; secondary 82A05 \endsubjclass

\endtopmatter

\bigskip
\subheading{ 1. Introduction} 
\numsec=1
\numtheo=1
\numfor=1

The rate of convergence to equilibrium is one of the main problems in
the theory of Markov processes. It has recently attracted the
attention of many authors in the context of symmetric conservative
particle systems in finite and infinite volume. In finite volume the
techniques used to obtain the rate of convergence to equilibrium rely
mostly on the estimation of the spectral gap of the generator. In
general, one shows that the generator of the particle system
restricted to a cube of length $N$ has a gap of order $N^{-2}$ in any
dimension. This estimate together with standard spectral arguments
permits to prove that the particle system restricted to a cube of size
$N$ decays to equilibrium in $L^2$ at the exponential rate
$\exp\{-ct/N^2\}$.  This approach has been successfully extended to
the infinite volume setting and permitted to prove $L^2$-polynomial
decay to equilibrium. The method, however, does not give any
information on the rate at which the system converges to equilibrium
starting from an arbitrary configuration or from arbitrary initial
distributions. In this article we consider the symmetric simple
exclusion process in which some explicit computations can be
performed. 

The symmetric simple exclusion process was introduced by Spitzer
(1970).  Informally one can describe the process following the so
called stirring representation. Fix a symmetric transition probability
$p$ on $\Bbb Z^d$~: $p(x,y) = p(y,x)$, $p(x,y)\ge 0$, $\sum_y p(x,y) =
1$ for all $x$ in $\Bbb Z^d$. We assume that the transition
probability is translation invariant, $p(x,y) = p(0,y-x)$ and that it
is indecomposable in the sense that for each $z$ in $\Bbb Z^d$, there
exists an integer $n$ and a sequence $0=z_0, \dots , z_n = z$ such
that $p(z_i, z_{i+1}) >0$ for all $i=0, \dots , n-1$. For each pair of
sites $x$, $y$ such that $p(x,y)>0$, consider a Poisson point process,
denoted by $N_{x,y} (t)$, with rate $p(x,y)$. Assume that these
processes are all independent. As initial state, fix some
configuration $\eta$ with at most one particle per site. Thus, the
configuration $\eta$ is a collection $\eta = \{\eta (x), \, x\in \Bbb
Z\}$, where $\eta (x) =1$ indicates that the site $x$ is occupied for
the configuration $\eta$ and $\eta (x)=0$ means that the site is
empty.  To obtain the state of the process from the initial states and
from the Poisson point processes, we proceed as follows. Each time the
Poisson process $N_{x,y}$ increases by $1$, we interchange the
variables $\eta(x)$ and $\eta (y)$. Notice that if both site $x$ and
$y$ are occupied or if both are vacant before the jump of $N_{x,y}$,
the configuration remains unchanged after the jump. In the other two
cases the modification can be interpreted as the jump of a particle
from the occupied site to the unoccupied site. Notice also that at
each later time each site is occupied by at most one particle and that
the total number of particles is conserved by the dynamics. This
explains why this stochastic dynamics is called the symmetric simple
exclusion process and in which sense it is conservative.

For $\rho\in [0,1]$, denote by $\nu_\rho$ the Bernoulli product measure with
density $\rho$. This is the probability measure on the configuration space
$\{0, 1\}^{\Bbb Z^d}$ obtained by putting a particle at each site with
probability $\rho$ independently of the other sites.  Liggett proved that
all the invariant measures for the symmetric simple exclusion process are
convex combinations of the Bernoulli measures $\nu_\rho$ (cf. Liggett (1985)).

When one tries to study convergence to equilibrium for arbitrary
initial configurations, the first issue is the very question: ``what
to prove?'' The first attempt is to try to prove the existence of a
function $h(t)$ that decreases to $0$ as $t\uparrow\infty$ and the
existence of a norm $|\! |\! | f |\! |\! |$ such that for all initial
configuration $\eta$ chosen according to the invariant measure
$\nu_\rho$ and for all cylinder function $f$
$$
|\delta_\eta S(t)f - \nu_\rho f |\,\le\,  \, |\! |\! | f |\! |\! | \, h(t) \; .
$$ Here $\{S(t), t\ge 0\}$ stands for the semigroup of the symmetric
simple exclusion process and $\delta_\eta$ for the probability measure
concentrated on $\eta$ so that $\delta_\eta S(t)$ is the
distribution of the process at time $t$ starting from
$\eta$. Also, for a bounded function $f$ and a probability measure
$\mu$, $\mu f$ stands for the expectation of $f$ with respect to
$\mu$.

This is clearly not true as for any fixed $\varepsilon>0$ and any
fixed $t$, one can always choose a set of configurations with
$\nu_\rho$ positive probability for which
$$ 
|\delta_\eta S(t) [\eta (0)] - \rho | \; \ge \; \varepsilon\; .
$$ 
Indeed, it is enough to consider configurations $\eta$ whose sites in a
cube of length $t$ around the origin are all occupied.  In general there is no
hope to have uniform almost sure convergence in conservative systems. In spin
flip systems the equilibrium is attained locally in an independent way for
distant regions. For this reason one can hope to get almost sure convergence
to equilibrium in spin-flip systems.

The second attempt is to fix a configuration $\eta$, to choose the density
$\rho$ depending on the initial configuration: $\rho^\eta_t = \delta_\eta
S(t) [\eta(0)]$ and to compute
$$
|\delta_\eta S(t)f - \nu_{\rho^\eta_t} f |\; .
\Eq(aa)
$$ 
This approach cannot give a bound better than $ c(f)t^{-1/2}$ in any
dimension. This is not satisfactory because, in view of the decay to
equilibrium in $L^2$, one expects to obtain estimates of order $t^{-d/2}$ in
dimension $d$. To check that this formulation can not give bounds better than
$t^{-1/2}$, consider the configuration $\eta$ such that $\eta (x)=1$ if and
only if $x_1\ge 1$. Here $x_1$ stands for the first coordinate of $x$. In this
case, a standard duality argument (that will be explained in section 2) gives
that $\rho^\eta_t = P[X^0_t\ge 1]$, where $X^0$ is a one dimensional symmetric
random walk that starts from the origin. Fix the cylinder function $f(\eta)
=\eta (e_1)$, where $e_i$, $1\le i\le d$, is the canonical basis of $\Bbb
R^d$. By duality and translation invariance, $\delta_\eta S(t)f = P[X^0_t\ge
0]$. In particular, for this cylinder function, \equ(aa) is equal to $P[X_t^0
=0]$, which is of order $t^{-1/2}$.

At this point we have two possibilities. We may of course impose some
regularity conditions on the initial configuration $\eta$ (assume for instance
that it is periodic) or to average the difference $|\delta_\eta S(t)f -
\nu_{\rho} f |$ (or a power of this difference) with respect to some measure
$\nu$ that has asymptotic density $\rho$ and some nice correlation
properties. The other possibility is to take advantage that we are in the
context of exclusion process, where all cylinder functions can be written as
linear combinations of functions of type $\prod_{x\in A} \eta (x)$ for some
finite set $A$. In this case, a natural quantity to investigate is 
$$
\delta_\eta S(t)  \prod_{x\in A} \big\{ \eta (x) - \rho^\eta_t(x)
\big\}\; ,
$$ 
where $\rho^\eta_t(x) = \delta_\eta S(t) \eta(x)$. These are the
so-called $v$-functions introduced by Ferrari, Presutti, Scacciatelli
and Vares (1991). A related quantity, also very natural in the context
of the symmetric simple exclusion process is the difference
$$
\delta_\eta S(t) \prod_{x\in A} \eta (x) \; -\; \prod_{x\in A} \delta_\eta
S(t) \eta (x)\; .
\Eq(bb)
$$ 

We shall follow the two directions just mentioned. We shall first prove that
the difference \equ(bb) can be expressed in terms of quantities related to
independent random walks with transition rate $p(x,y)$. This is the
content of Theorem \equ(basic). This fact together with
some elementary bounds on the transition probability of symmetric random
walks will permit to obtain sharp estimates on the integral of powers of
\equ(aa) with respect to translation-invariant measures that have density
$\rho$ and polynomial decaying correlations. This is the content of Theorem
\equ(theo).

The comparison between $n$ random walks interacting by exclusion and $n$
independent random walks was studied by Bertein and Galves (1977), De Masi,
Ianiro and Presutti (1982), De Masi and Presutti (1983), De Masi, Ianiro,
Pellegrinotti and Presutti (1984), Ferrari and Goldstein (1988), Ferrari,
Presutti, Scacciatelli and Vares (1991) and Andjel (1994). We give a short and
unified presentation which includes most of the above results.

A simple exclusion process in a finite box can be understood as a
random walk in a finite set.  Convergence to equilibrium for
finite-state Markov processes have received new attention lately. We
quote for instance Aldous (1983), Diaconis (1988), Diaconis and
Stroock (1991). For the symmetric simple exclusion in a finite box,
Quastel (1992) computed the spectral gap of the generator, which gives
the $L^2$ exponential rate of convergence to equilibrium in finite
volume. Using this result, Bertini and Zegarlinski (1998) proved
polynomial $L_2$-convergence to equilibrium (time-correlation decay
for the system in equilibrium) in infinite volume.  Janvresse, Landim,
Quastel and Yau (1999) proved a analogous result for the symmetric
zero range process. 

Cancrini and Galves (1995) obtained an upper bound
for the rate of convergence of symmetric simple exclusion processes
starting either from a periodic configuration or from a stationary
measure satisfying mixing conditions. This result was extended to the
one-dimensional nearest-neighbor zero-range process with rate $g(k) =
\bold 1\{k\ge 1\}$ by Galves and Guiol (1997).

\subheading{ 2. Notation and results}
\numsec=2
\numtheo=1
\numfor=1

Let $p(x,y)$ be a symmetric, translation-invariant, irreducible, 
non-negative real matrix on $\Bbb Z^d$ such that
$$
\sum_y |y|^2 p(0,y)<\infty\; .
$$
To each pair of sites $(x,y)$ of $\Z^d$ attach a Poisson process of
rate $p(x,y)$ denoted by $N_{x,y}(t)$. We adopt the following
convention: $dN_{x,y}(t)=1$ when there is an event of the
corresponding Poisson process at time $t$. As a function of these
Poisson processes, for $s<t$, we define a family of bijections
of $\Z^d$ in itself, $\xi^s_t\colon \Z^d \to \Z^d$,  such that
$$
\xi^s_s(x) = x \quad \text{and}\quad 
d\xi^s_t = \sum_{x,y} \{ (\xi^s_t)^{x,y} - \xi^s_t \} dN_{x,y}^s(t) \; ,
$$
where for any bijection $\xi\colon \Z^d\to\Z^d$,
$$
\xi^{x,y}(z) =\cases \xi(z) &\text{if } z\neq x,y\, ,\cr
\xi(y) & \text{if } z=x\, ,\cr
\xi(x) & \text{if } z=y\, \cr
\endcases
$$
and $N_{x,y}^s(t) = N_{x,y}(s+t) - N_{x,y}(s)$.
$\xi_\cdot^s$ is thus the identity at time $s$ and if there is a mark
of the Poisson process $N_{x,y}$ at time $s+t$, $\xi^s (x)$ and
$\xi^s(y)$ change their value in the following way~: $\xi^s_t(x) =
\xi^s_{t-} (y)$, $\xi^s_t(y) = \xi^s_{t-} (x)$. In particular, for any
fixed $x$, the process $\xi^s_\cdot (x)$ is a stochastic process that
may have very long jumps even if the rates $p(\cdot, \cdot)$ are
nearest-neighbor. 

The simple exclusion process with initial configuration $\eta$ is
defined as
$$
\eta^\eta_t(x) = \eta(\xi^0_t (x))
$$
for all $t\ge 0$.
Denote $S(t)f(\eta)= \E f(\eta^\eta_t)$. Here and below $\P$ and
$\E$ stand for the probability and the expectation with respect to the
Poisson point processes. The set of extremal invariant measures for
this process is the set of Bernoulli product measures $\nu_\rho$
indexed by the density $\rho\in [0,1]$.
The main results of this article concern the convergence to those measures from
different initial conditions. They are stated in the next theorem.

Fix $\alpha>0$. We say that a measure $\nu$ on $\{0,1\}^{\Z^d}$ has
$\alpha$-decaying correlations if there exist a finite constant $C$
such that 
$$
\Big |\, E_\nu[fg] - E_\nu[f] E_\nu[g] \, \Big| \; \le \; C \Vert f
\Vert_\infty \Vert g \Vert_\infty \frac 1 {d(\Lambda_f ,
  \Lambda_g)^\alpha}
$$
for every cylinder functions $f$, $g$. In this formula,
$\Lambda_f$, $\Lambda_g$ stands for the support of the cylinder
functions $f$, $g$ and $d(A,B)$ for the distance between two subsets
$A$, $B$ of $\Z^d$.

\proclaim{\Theorem(theo)} Let $\nu$ be a translation-invariant
probability measure in $\{0,1\}^{\Z^d}$ with
density $\rho$ (that is $E_\nu [\eta(x)] = \rho$) and with 
$\alpha$-decaying correlations for some $\alpha > d$. Then, 
for each cylinder function
$f$ there exists a constant $c(f)$ such that for all $t>0$:
\item{(a)} Weak convergence:
$$
|\nu S(t) f - \nu_\rho f| \,\le\, c(f)\, (1+t)^{-d/2}
$$

\item{(b)} $L_p$ convergence for $1\le p\le 2$: 
$$
\int \nu(d\eta)\,|\delta_\eta S(t) f - \nu_\rho f|^p \, \le\, c(f)\,
(1+t)^{-dp/4} \, .
$$
\item{(c)} $L_p$ convergence for $p\ge 2$: 
$$
\int \nu(d\eta)\,|\delta_\eta S(t) f - \nu_\rho f|^p\, \le\, c(f)\, 
(1+t)^{-d/2} \; .
$$
\endproclaim 

\noindent{\bf Remarks.} Part (a) in the theorem above improves the
bound obtained by Cancrini and Galves (1995) in two directions: on the
one hand, we only ask for polynomial decay of correlations in the
initial measure, while [CG] requires exponential decay. On the other
hand, in [CG] the upper bound includes a (spurious) $\log t$ in the
numerator.  Part (b) generalizes a result obtained by Bertini and
Zegarlinski (1996) in two senses: [BZ] proved (b) for $p=2$ with
product initial measure, while part (b) requires only the initial
measure to have polynomial decay of correlations, which is of course
satisfied by $\nu_\rho$, and part (b) is proved for $1\le p\le 2$.
Finally, since $\nu_\rho$ is an equilibrium measure, part (b) can be
read as the time-decay of correlations of the system in equilibrium.
%It follows from the proof of this theorem that for $p\ge 2$,
%$$
%\int \nu(d\eta)\,|\delta_\eta S(t) f - \nu_\rho f|^p\, \le\, c(f)\, 
%(1+t)^{-d/2} \; .
%$$

The proof of Theorem \equ(theo) is based on the self-duality of the
symmetric simple exclusion process that we now explain.  For $0\le
s\le t$, let $X^{x,t}_s$ denote the position at time $t-s$ of the
particle sitting at $x$ at time $t$~: $X^{x,t}_s = \xi^{t-s}_s(x)$. In
contrast with $\xi^s_\cdot (x)$, $X^{x,t}_\cdot $ is a random walk with
transition rate $p(\cdot , \cdot)$. Moreover, for any subset $A$
of $\Bbb Z^d$, $\{X^{x,t}_\cdot ,\, x\in A\}$ evolves as symmetric
exclusion random walks and
$$
\P[\eta^\eta_t (x) = 1, \, x\in A] \; =\; \P[\eta (X^{x,t}_t) = 1 ,\,
x\in A]
$$
for every finite subset $A$.  This is the so-called self-duality
relation of the symmetric simple exclusion process. 

In the following theorem we compare the evolution of $n$ particles
interacting by exclusion with $n$ independent particles. The proof is
inspired in a similar result by Ferrari and Goldstein (1988), where
exclusion processes with birth and deaths were considered, and a
result by Ferrari, Presutti, Scacciatelli and Vares (1991). This can
be seen as a probabilistic version of the ``integration by parts
formula''. See Proposition 8.1.7 of Liggett (1985) and display (1) of
Andjel (1994), for instance.

\proclaim{\Theorem(basic)} For  any vector $\underline x =
(x_1,\dots,x_n)$,
$$
\aligned
\E\, &\prod_{i=1}^n \eta^\eta_t(x_i) - \prod_{i=1}^n
  \E\, \eta^\eta_t(x_i) \cr
&= - \int_0^t \sum_{\underline y} \P(\underline
  X^{\underline x,t}_s = \underline y)\, \sum_{i<j} p(y_i,y_j)\,
    (\rho^\eta_{t-s}(y_i)-\rho^\eta_{t-s}(y_j))^2 \,\prod_{k\neq i,j}
    \rho^\eta_{t-s}(y_k)
\endaligned
$$
for all $t$ and any configuration $\eta$. In the above formula
the summation is carried over all $\underline y = (y_1, \dots, y_n)$
such that $y_i\neq y_j$ for $i\neq j$ and
$$
\rho^\eta_t(y) = \E\, \eta^\eta_t(y).
$$
\endproclaim

In the next theorem we apply the previous result.

\proclaim{\Theorem(11)} Assume $p(\cdot,\cdot)$ has a finite second
moment. 
{\parindent .7cm
\item{(a)} There exists a finite constant $C$ such that for all $t\ge 0$,
$$
\sup_\eta \Big|\E\,\prod_{i=1}^n \eta^\eta_t(x_i) - \prod_{i=1}^n
  \E\,\eta^\eta_t(x_i)\Big| \le  C \, R_d (t)\; ,
$$
where $R_1(t) = \log (1+t)/\sqrt{1+t}$, $R_2 (t)= \log (1+t)/(1+t)$
and $R_d (t) = (1+t)^{-1}$ for $d\ge 3$.
\item{(b)} For each positive function $\varphi:\R_+\to\R_+$
let  
$$
\bold X_\varphi := \Big\{ \eta: \text{ for all }x,y\in\Z^d,\;
|x-y|^{-2}(\rho^\eta_t(x) - \rho^\eta_t(y))^2 \;\le\; \varphi(t)
\Big\}\; . 
$$ 
Then, for each decreasing function $\varphi$ such that $\varphi (t)
\le C_0 (1+t)^{-(d+2)/2}$ for some finite constant $C_0$, there exists a
finite constant $C_1$ such that for all $t\ge 0$
$$
\sup_{\eta \in\bold X_\varphi} \Big|\E\,\prod_{i=1}^n
\eta^\eta_t(x_i) - \prod_{i=1}^n \E\,\eta^\eta_t(x_i)\Big| \le C_1 {1\over
  (1+t)^{d/2}}\; ;
$$
\item{(c)} Let $\nu$ satisfy the conditions of Theorem
\equ(theo). Then there exists a finite constant $C$ such that
for all $t\ge 0$
$$
\nu\Big|\E\,\prod_{i=1}^n \eta^\eta_t(x_i) - \prod_{i=1}^n
  \E\,\eta^\eta_t(x_i)\Big| \le  C { 1 \over (1+t)^{d/2}}\; .
$$
}
\endproclaim

\noindent{\bf Remarks. } Part (a) was obtained by Ferrari, Presutti,
Scacciatelli and Vares (1988). Let $\eta$ be a periodic configuration on $\Bbb
Z^d$. It is easy to show that $\eta$ belongs to $\bold X_\varphi$ for some
$\varphi$ satisfying the assumptions of Theorem \equ(11) (b). Part (b)
includes a result in Landim (1999), where the case of finite initial $\eta$ is
considered. 

\subheading{ 3. Proofs}
\numsec=3
\numtheo=1
\numfor=1
 
Before proving the theorems, we state an estimate that will be needed
several times in the sequel. It is based on the classical negative
correlations property of the symmetric exclusion process.

\proclaim{\Lemma(bes)}
Let $\Phi:\Z^d\setminus\{0\}\to\R^+$ be a summable function. Then,
there exists a 
constant $C$ such that for all set of
different sites $\{x_1,\dots, x_n\}$ and for all $t>0$,  
$$
\sum_{j\neq k} \sum_{y_1, \dots , y_n} \Bbb P \big[ \{
X^{x_1,t}_t, \dots ,X^{x_n,t}_t\} = \{y_1, \dots ,y_n\} \big] \Phi(y_j
- y_k) \; \le\; C \,(1+t)^{-d/2},\Eq(bes1)
$$
where the second sum is carried over the set $\{y_1, \dots , y_n: y_a\neq
y_b,$ for $a\neq b\}$.
\endproclaim 

\demo{Proof} If we fix $j\neq k$
and sum over $y_i$ for $i\neq j$, $k$, we obtain that the left hand
side of \equ(bes1) is equal to
$$
\sum_{j\neq k} \sum_{ y_j\neq y_k} \Bbb P \big[ \{ X^{x_1,t}_t,
\dots ,X^{x_n,t}_t\} \supset \{y_j, y_k\} \big] \Phi(y_j - y_k) \; \cdot
\Eq(caa)
$$ 
The previous probability can be decomposed as
$$
\sum_{a\neq b} \Bbb P \big[ \{ X^{x_a,t}_t ,X^{x_b,t}_t\} = \{y_j,
y_k\} \big].
$$
By the correlation inequality between symmetric exclusion random
walks and symmetric independent random walks (cf. [L], Theorem
VIII.1.7), each of theses probabilities is bounded above by
$$
\aligned
&  \Bbb P \big[ X^{x_a,t}_t \in \{ y_j,  y_k\} \big]  \Bbb P \big[
X^{x_b,t}_t \in \{ y_j,  y_k\} \big]  \\
&\qquad = \; \big\{ p_t (y_j - x_a) + p_t(y_k - x_a) \big\} \big\{ p_t
(y_j - x_b) + p_t(y_k - x_b) \big\}  
\endaligned
$$ 
provided $p_t(x,y)$ stands for the probability of a continuous-time
random walk with transition probability $p(\cdot , \cdot)$ to be at
$y$ at time $t$ if it starts from $x$ at time $0$ and provided
$p_t(x)$ stands for $p_t(0,x)$. The sum \equ(caa) is thus bounded
above by
$$
\sum_{a\neq b} \sum_{j\neq k} \sum_{y_j\neq y_k} \Phi(y_j - y_k)\, 
\big\{ p_t (y_j - x_a) + p_t(y_k - x_a) \big\} 
\big\{ p_t(y_j - x_b) + p_t(y_k - x_b) \big\} \; .
$$
Since $\Phi$ is summable and since $p_t(x)$ is bounded above by
$C/(1+t)^{d/2}$ for some finite constant that depends only on
$p(\cdot)$, this expression is bounded above by $C(n,p) (1+t)^{-d/2}$,
which concludes the proof of the lemma.\ \qed
\enddemo

\noindent{\bf Proof of Theorem \equ(theo) (a)} Let $\nu$ be a
translation-invariant probability measure that has $\alpha$-decaying
correlations.  By duality, for any $n\ge 1$ and any distinct sites
$x_1, \dots, x_n$,
$$
\aligned 
& \Big\vert \Bbb P_\nu \big[ \eta _t (x_i) = 1\, , \; i=1,
\dots , n \big] - \rho^n \Big\vert \\
& =\; \Big\vert \sum_{\scriptstyle y_1, \dots , y_n\atop
  \scriptstyle y_i \neq y_\ell } \Bbb P \big[ \{ X^{x_1,t}_t, \dots
,X^{x_n,t}_t\} = \{y_1, \dots ,y_n\} \big] \Big\{ \nu \{ \eta (y_i) =1 ,
i=1, \dots , n \} - \rho^n \Big\} \Big\vert \; .  
\endaligned
$$
Since $\nu$ is translation-invariant with mean $\rho$ and has
$\alpha$-decaying correlations, there exists a finite
constant $C$ such that
$$
\Big\vert \nu \{ \eta (y_i) =1 , i=1, \dots
, n \} - \rho^n \Big\vert \; \le \; C \sum_{j\neq k} \frac 1 {|y_j -
  y_k|^\alpha}\; \cdot
$$
The right hand side of the previous expression is thus bounded
above by
$$
C \sum_{j\neq k} \sum_{\scriptstyle y_1, \dots , y_n\atop
  \scriptstyle y_i \neq y_\ell} \Bbb P \big[ \{ X^{x_1,t}_t, \dots
,X^{x_n,t}_t\} = \{y_1, \dots ,y_n\} \big] \frac 1 {|y_j - y_k|^\alpha}
\; \cdot
$$
To conclude the proof it remains to apply Lemma \equ(bes) for the
function $\Phi(y) = |y|^{-\alpha}$. We can do this because in the
above expression $y_j\neq y_k$.
  \qed 
\medskip

\noindent{\bf Proof of Theorem \equ(basic).}
Fix $n$ distinct sites $x_1,\dots,x_n$ and define
$\{Y^{x_i,s}_{i,t}:i=1,\dots,n\}$ as a family of independent random
processes with the same marginal distribution as $X^{x_i,s}_{t}$
respectively.

We realize the motion of the $Y$ process by considering a family of
Poisson marks as the one defined by $N_{x,y}(t)$, called
$N^0_{x,y}(t)$, independent of the precedent one as follows.
Processes $X$ uses only the $N$ marks. Process $Y$ uses the $N$ marks
if only one $Y$ particle is concerned by the jump. This means that if
there is a $Y$ particle at $x$ and no $Y$ particle at $y$ at time $t$,
then the process uses the marks of $N_{x,y}(t)$ and ignores the marks
of $N_{x,y}^0(t)$. If the jump concerns two $Y$ particles, say
particles $i$ and $j$ with $i<j$, then the $Y_i$ particle uses the $N$
marks and the $Y_j$ particle uses the $N^0$ marks to jump over the
position of the $Y_i$ particle and the $N$ marks to jump to any other
position. Thus if particle $Y_i$ is at $x$ and particle $Y_j$ at time
$t$ with $i<j$, particle $Y_i$ uses the marks of $N_{x,y}(t)$ while
particle $Y_j$ uses the marks of $N_{x,y}^0(t)$.

We just gave a {\sl coupling} between a system of $n$ exclusion and
$n$ independent particles. This means that we realized the
two processes in the same probability space (the one generated by the
product of the Poisson processes $\{N_{x,y}(t), \, N_{x,y}^0(t), \,
;\; x, y \in \Z^d,\, t\ge 0\}$) in such a way that the marginal
distributions are those desired for both processes. We continue using
$\P$ and $\E\,$ for the probability and expectation with respect to the
product of the Poisson processes.

By definition of the symmetric simple exclusion process,
$$
\E\, \Big[ \prod_{i=1}^n \eta^\eta_t(x_i ) \Big] - \prod_{i=1}^n \E\,
\Big[\eta^\eta_t(x_i) \Big] \; =\; \E\, \Big[ \prod_{i=1}^n
\eta(X^{x_i,t}_{t}) - \prod_{i=1}^n \eta(Y^{x_i,t}_{i,t})\Big]\; .
$$

Let $T_1$ be the first instant that two $Y$ particles occupy the same
site. We say that a {\sl collision} occurred at that time.  Before the
collision each $X$ particle occupies the same place of the
corresponding $Y$ particle.  For $s<t$ let
$$
I(x_1,\dots,x_n,s,t)=\prod_{i=1}^n \eta(X^{x_i,t}_{s}) 
- \prod_{i=1}^n \eta(Y^{x_i,t}_{i,s})\; .
$$
We want to compute the expectation of $ I(x_1,\dots,x_n,t)
= I(x_1,\dots,x_n,t,t)$. If $t<T_1$, $I(x_1,\dots,x_n,t)$ is zero 
because the trajectories of the $X$ and $Y$ process coincide.
On the other hand, since for all $0\le s\le t$,
$Y^{x_i,t}_{i,t} = Y^{Y^{x_i,t}_{i,s}, t-s}_{i,t-s}$,
$X^{x_i,t}_{t} =X^{X^{x_i,t}_{s}, t-s}_{t-s}$, on the
set $\{T_1\le t\}$,
$$
I(x_1,\dots,x_n,t) \; =\; \prod_{i=1}^n \eta(X^{X^{x_i,t}_{T_1}, t-T_1}
_{t-T_1}) - \prod_{i=1}^n \eta(Y^{Y^{x_i,t}_{i,T_1},
t-T_1}_{i,t-T_1})\; .
\Eq(fl0)
$$

Let
$$
E_1 = 
\cases 
1 &\text{if the collision at $T_1$ occurs due to a $N$ mark}\cr
0 &\text{if the collision at $T_1$ occurs due to a $N^0$ mark}
\endcases
$$
and let $Z_1\in \{(i,j),\; 1\le i<j\le n\}$ stand for the labels of
the particles involved in the collision at time $T_1$~:
$$
Z_1 = (i,j) \quad \text{if and only if} \quad Y^{x_i,t}_{i,T_1}=
Y^{x_j,t}_{j,T_1} \quad \text{and}\quad i<j\; .
$$
Assume that the first collision is
due to a collision between particles $i$ and $j$ with $i<j$. In this
case, at time $T_1$, the position of particles $Y^{x_k,t}_k$ and
$X^{x_k,t}$ coincide for $k\neq i$, $j$. Moreover, if the collision 
occurred due to a $N$
mark, i.e., due to a jump of particle $Y_i$ over particle $Y_j$,
$Y^{x_i,t}_{i,T_1} = X^{x_j,t}_{j,T_1} = Y^{x_j,t}_{j,T_1}$. In the
case where the collision occurred due to a $N^0$ mark, i.e., due to a
jump of the $Y_j$ particle over $Y_i$, a similar identity holds with
the roles of $i$ and $j$ exchanged.  In particular, on the set
$\{T_1\le t\}$, 
$$
\aligned
& \prod_{k=1}^n \eta(Y^{Y^{x_k,t}_{T_1},t-T_1}_{k,t-T_1}) \; =\; 
\sum_{i<j}\one\{Z_1=(i,j)\}
\prod_{k\notin\{i,j\}} \eta(Y^{X^{x_k,t}_{T_1},t-T_1}_{k,t-T_1}) \;
\times \cr
& \qquad\qquad\qquad\qquad\qquad\qquad \times \Bigl\{ \one\{E_1=0\}
  \eta(Y^{X^{x_i,t}_{T_1},t-T_1}_{i,t-T_1} )
\eta(Y^{X^{x_i,t}_{T_1},t-T_1}_{j,t-T_1}) \cr 
&\qquad\qquad\qquad\qquad\qquad\qquad\qquad + \; \one\{E_1=1\}
  \eta(Y^{X^{x_j,t}_{T_1},t-T_1}_{i,t-T_1} )
\eta(Y^{X^{x_j,t}_{T_1},t-T_1}_{j,t-T_1}) 
\Bigr\} \; .
\endaligned
$$
We now add and subtract in the right hand side
$\eta(Y^{X^{x_i,t}_{T_1},t-T_1}_{i,t-T_1} )
\eta(Y^{X^{x_j,t}_{T_1},t-T_1}_{j,t-T_1})$ to recover
$\prod_{1\le k\le n} \eta(Y^{X^{x_k,t}_{T_1},t-T_1}_{k,t-T_1})$.
After this step we obtain that on the set
$\{T_1\le t\}$, 
$$
\aligned
& \prod_{k=1}^n \eta(Y^{Y^{x_k,t}_{T_1},t-T_1}_{k,t-T_1}) \\
& \quad  = \; \prod_{k=1}^{n} \eta(Y^{X^{x_k,t}_{T_1},t-T_1}_{k,t-T_1})
\; +\;\sum_{i<j}\one\{Z_1=(i,j)\}
\prod_{k\notin\{i,j\}} \eta(Y^{X^{x_k,t}_{T_1},t-T_1}_{k,t-T_1}) \;
\times \cr
& \quad \times \Bigl\{ \one\{E_1=0\}
 \Big[  \eta(Y^{X^{x_i,t}_{T_1},t-T_1}_{i,t-T_1} )
\eta(Y^{X^{x_i,t}_{T_1},t-T_1}_{j,t-T_1}) -
\eta(Y^{X^{x_i,t}_{T_1},t-T_1}_{i,t-T_1} )
\eta(Y^{X^{x_j,t}_{T_1},t-T_1}_{j,t-T_1}) \Big] \cr 
&\qquad + \; \one\{E_1=1\}
\Big[  \eta(Y^{X^{x_j,t}_{T_1},t-T_1}_{i,t-T_1} )
\eta(Y^{X^{x_j,t}_{T_1},t-T_1}_{j,t-T_1})  -
\eta(Y^{X^{x_i,t}_{T_1},t-T_1}_{i,t-T_1} )
\eta(Y^{X^{x_j,t}_{T_1},t-T_1}_{j,t-T_1}) \Big]  
\Bigr\} \; .
\endaligned
$$
Denote the right hand side of this expression by $g(E_1, T_1, Z_1,
\underline X^{\underline x, t}_{T_1}, \underline Y^
{\underline X^{\underline x, t}_{T_1}, t -T_1}_{t-T_1})$. It
follows from this identity and \equ(fl0) that
$$
I(\underline x, t) \; =\; \bold 1\{T_1\le t\}
I(\underline X^{\underline x, t}_{T_1}, t-T_1) \; - \;
\bold 1\{T_1\le t\} g(E_1, T_1, Z_1,
\underline X^{\underline x, t}_{T_1}, \underline Y^
{\underline X^{\underline x, t}_{T_1}, t -T_1}_{t-T_1})\; .
$$

For a positive integer $\ell$, define ${\Cal G}_\ell=
\sigma\{T_m,Z_m,X^{x_1,t}_{T_m}, \dots, X^{x_n,t}_{T_m}:
m=1,\dots,\ell\}$.  Since $E_1$ is a Bernoulli random variable with
parameter $1/2$ independent of ${\Cal G}_1$ and since the $Y$
particles evolve independently, by the strong Markov property on the
set $\one\{T_1\le t\}$ we get that
$$
\E\,( g(E_1, T_1, Z_1,
\underline X^{\underline x, t}_{T_1}, \underline Y^
{\underline X^{\underline x, t}_{T_1}, t -T_1}_{t-T_1}) | {\Cal G}_1)\;
= \; g(T_1,Z_1,X^{x_1,t}_{T_1},\dots,X^{x_n,t}_{T_1})\; ,
$$
where
$$
g(u,\{i,j\},x_1,\dots,x_n) = {1\over 2}
\,\left[\E\,\left(\eta(Y^{x_i,t}_{i,t-u})\right)-
  \E\,\big(\eta(Y^{x_j,t}_{j,t-u})\big)\right]^2 
\prod_{k\neq i,j} \E\,\left(\eta(Y^{x_k,t}_{k,t-u})\right)\; .
$$
 From the two previous identities, we obtain that
$$
\aligned
\E\,( I(x_1,\dots ,x_n, t) | {\Cal G}_1) \; & =\;
\one\{T_1\le t\} I(X^{x_1,t}_{T_1},\dots,X^{x_n,t}_{T_1}, t-T_1) \\
& \quad-\;\one\{T_1\le t\} g(T_1,Z_1,X^{x_1,t}_{T_1},\dots,
X^{x_n,t}_{T_1})\; .
\endaligned
$$
Repeating the argument for $\ell \ge 2$, we get the following
expression for the expectation of~$I$:
$$
\aligned
\E\,(I&(x_1,\dots,x_n,t)) = - \E\, \sum_{\ell= 1}^{M(t)} 
g\big(T_\ell,Z_\ell,X^{x_1,t}_{T_\ell},\dots,X^{x_n,t}_{T_\ell}\big)
\; ,
\endaligned
$$ 
where $T_\ell, E_\ell, Z_\ell$ are defined as $T_1$, $E_1$, $Z_1$
inductively and $M(t)= \sum_{\ell\ge 1} \one\{T_\ell\le t\}$ is the
number of collisions occurred by time $t$.  Notice that since $g$ is
positive, the above expression implies immediately that the
distribution of $\eta_t$ is a measure with negative correlations for
any initial $\eta$.

For $i<j$, denote by $T^{i,j}_\ell$ the instant of the $\ell$-th collision 
of particles $i,j$ and by $M^{i,j}(t)$ the number of collisions up to time
$t$ of particles $i,j$ so that $M(t) = \sum_{i<j} M^{i,j}(t)$ and
$$
 \E\,(I(x_1,\dots,x_n,t))\;=\; - \sum_{i<j} \E\, \sum_{\ell=
  1}^{M^{i,j}(t)} g^{i,j}\big(T^{i,j}_\ell,
X^{x_1,t}_{T^{i,j}_\ell},\dots,X^{x_n,t}_{T^{i,j}_\ell} \big) \; ,
\Eq(10)
$$
where 
$$
g^{i,j}(u,\underline y) \;=\; {1\over 2}
\,\left[\E\,\left(\eta(Y^{y_i,t}_{i,t-u})\right)-
  \E\,\big(\eta(Y^{y_j,t}_{j,t-u})\big)\right]^2 \,\prod_{k\neq i,j}
\E\,\eta(Y^{y_k,t}_{k,t-u}) \; .
\Eq(09)
$$

For $i<j$, $M^{i,j}(t)$ is a Poisson process with rate
$2 \sum_{\underline y} \one \{\underline X_t = \underline y\}
p(y_i,y_j)$. Moreover, for any function $h (u, \underline x)$,
$$
\sum_{\ell=1}^{M^{i,j}(t)} h \big(T^{i,j}_\ell,
\underline X^{\underline x,t}_{T^{i,j}_\ell} \big)
\; =\; \int_0^t h \big( s ,\underline X^{\underline x,t}_{s} \big)
\, d M^{i,j}(s)\; .
$$
In particular, taking expectations,
$$
\E\, \sum_{\ell= 1}^{M^{i,j}(t)}
 g^{i,j}\big(T^{i,j}_\ell,
\underline X^{\underline x,t}_{T^{i,j}_\ell}\big) \;=\; 2 \int_0^t
\sum_{\underline y} p(y_i,y_j) \P(\underline X^{\underline
  x}_s=\underline y)\,g^{i,j}(s,\underline y)\, ds \; .
$$
This together with \equ(10), \equ(09) concludes the proof of the theorem.
\qed 

\bigskip
\noindent{\bf Proof of Theorem \equ(11).}
To show item (a) we use Theorem \equ(basic). Since we are taking supremum
over $\eta$, we can cancel the product $\prod_{k\ne
  i,j}\rho^\eta_{t-s}(y_k)$. We get
the following upper bound 
$$
\aligned
\sum_{i<j}\int_0^t \sum_{y_1,\dots, y_n} p(y_i,y_j)\,
\P(\{  X^{x_1,t}_s
,\dots,X^{x_n,t}_s \}&= \{y_1,\dots,y_n\})\, \times \cr
&\times    (\rho^\eta_{t-s}(y_i)-\rho^\eta_{t-s}(y_j))^2 \,ds 
\endaligned
\Eq(12)
$$
Denote by $p_t(x,y)$ the probability for a random walk with transition
probability $p(\cdot , \cdot)$ starting at $x$ to be at $y$ a time $t$. 
By Lemma \equ(a) below, there exists a universal constant $C_1$ such that
$$
|\rho^\eta_{s}(y_i)-\rho^\eta_{s}(y_j)|\;=\;\Big|\sum_{z\in \Bbb Z^d}\eta(z)
\,(p_s(y_i,z)-p_s(y_j,z))\Big| \;\le\; C_1\,{|y_i-y_j|\over \sqrt {1+s}}
$$
for all $s\ge 0$. By Lemma \equ(bes) with $\Phi(y)= p(0,y) |y|^2$,
which is summable by hypothesis, \equ(12) is bounded above by
$$
C \int_0^t (1+s)^{-d/2 }\,{1\over 1+t-s}\,ds\;\le\;
C \, R_d (t) \,.
$$ 
Here and below $C$ stands for a finite constant that may change
from line to line.  This shows item (a).  To show (b), again we can
cancel the product $\prod_{k\ne i,j}\rho^\eta_{t-s}(y_k)$ and need to
bound \equ(12) for $\eta$ in $\bold X_\varphi$. For those $\eta$'s we
have
$$
(\rho^\eta_{t-s}(y_i)-\rho^\eta_{t-s}(y_j))^2\;\le\;
\varphi(t-s)\,|y_i-y_j|^2. 
$$ 
Applying Lemma \equ(bes) with $\Phi(y)= p(0,y) |y|^2$, we bound
\equ(12) by
$$
C\int_0^t (1+s)^{-d/2}\,\varphi(t-s)\,ds\;\le\; C (1+t)^{-d/2}
$$
because $\varphi (t)$ was assumed to be bounded above by
$C(1+t)^{-(d+2)/2}$. This proves (b).

To show (c), as before, we need to compute the expectation with
respect to $\nu$ of \equ(12).
Hence using Fubini  we need to compute $\nu
(\rho^\eta_{t-s}(y_i)-\rho^\eta_{t-s}(y_j))^2$. Since $\rho^\eta_t(y)
= \sum_z \eta(z) p_t(y,z)$, developing the square we get 
$$
\aligned
\nu(&\rho^\eta_{t-s}(y_i)-\rho^\eta_{t-s}(y_j))^2\cr
&= \sum_{z_1,z_2}
\nu(\eta(z_1)\eta(z_2))\,(p_t(y_i,z_1)-p_t(y_j,z_1))\,
(p_t(y_i,z_2)-p_t(y_j,z_2))\cr
& =\rho\sum_z(p_t(y_i,z)-p_t(y_j,z))^2 \cr
&\qquad+ \sum_{z_1\ne z_2}
(\nu\eta(z_1)\eta(z_2)- \rho^2) \,(p_t(y_i,z_1)-p_t(y_j,z_1))\,
(p_t(y_i,z_2)-p_t(y_j,z_2))\cr
&\qquad+ \sum_{z_1\ne z_2}
\rho^2\,(p_t(y_i,z_1)-p_t(y_j,z_1))\,
(p_t(y_i,z_2)-p_t(y_j,z_2)) \; .
\endaligned
$$ 
It follows form Lemma \equ(a) below that the sum of the first and
third line of the previous expression is bounded above by $C \rho
(1-\rho) |y_i-y_j|^2 (1+t)^{-(d+1)/2}$, while the second, by Schwarz
inequality, is bounded above by 
$$
\sum_{z_1\ne z_2} | \nu\eta(z_1)\eta(z_2)- \rho^2|
\,\{p_t(y_i,z_1)-p_t(y_j,z_1)\}^2\, .
$$
Since $\nu$ has $\alpha$-decaying correlations, by Lemma \equ(a), this
expression is bounded above by $C |y_i-y_j|^2 (1+t)^{-(d+1)/2}$.
Therefore, \equ(12) is less than or equal to
$$
C(\rho) \sum_{i<j}\int_0^t \sum_{y_1,\dots, y_n} p(y_i,y_j)
|y_i-y_j|^2 \, \P( \underline X^{\underline x,t}_s = \underline y)
\frac 1{(1+t-s)^{(d+1)/2}} \,ds \; .
$$
It remains to apply Lemma \equ(bes) to the
function $\Phi(y)= p(0,y) |y|^2$ to obtain that the previous
expression is bounded above by
$$
C(\rho) \int_0^t \frac 1{(1+s)^{d/2}}
\frac 1{(1+t-s)^{(d+1)/2}} \,ds \; .
$$ 
This concludes the proof of the theorem in dimension $d\ge 2$. In
dimension $1$, we need \equ(b) instead of Lemma \equ(a) to prove the
estimate. \qed \medskip

We turn now to Theorem \equ(theo).

\medskip \noindent{\bf Proof of Theorem \equ(theo).} We already proved
part (a). In order to prove (b), fix $1\le p\le 2$. It is enough to
consider cylinder functions of type $f_A= \prod_{x\in A} \eta (x)$ for
a finite subset $A$ of $\Bbb Z^d$.  For such a cylinder function,
$\delta_\eta S(t) f_A = \E\, \prod_{x\in A} \eta^\eta_t (x)$ and
$\nu_\rho f_A = \rho^{|A|}$, provided $|A|$ stands for the number of
sites of $A$. Since $0\le f_A \le 1$ and since $a^p \le a$ for
$0\le a\le 1$, by Theorem \equ(11) (c),
$$
\int \nu(d\eta)\,|\delta_\eta S(t) f_A - \nu_\rho f_A|^p \; \le \;
2^{p-1} \int \nu(d\eta)\, \Big | \prod_{x\in A} \E\,\eta^\eta_t(x)  
- \rho^{|A|} \, \Big | ^p \; + \; C (1+t)^{-d/2} 
$$
for some finite constant $C$. Introducing intermediary terms, the
first term on the right hand side is bounded above by
$$
C(A) \sum_{x\in A} \int \nu(d\eta)\, \Big| \E\,\eta^\eta_t(x)  
- \rho \, \Big| ^p\; .
$$ 
Since $p\le 2$, By H\"older inequality, the previous expression is
bounded above by
$$
C(A) \sum_{x\in A} \Big\{ \int \nu(d\eta)\, \Big( \E\,\eta^\eta_t(x)  
- \rho \Big) ^2 \Big\}^{p/2}\; .
$$ To show that this expression is bounded above by $C(A)
(1+t)^{-dp/4}$ it is enough to expand the square and to recall that
$\nu$ has $\alpha$-decaying correlations and that $p_t(x)$ is bounded
above by $C t^{-d/2}$ uniformly in $x$. The proof of (c) is exactly
the same. This concludes the proof of the theorem. \qed \medskip

We conclude this section with an estimate on the transition
probability of symmetric random walks.

\proclaim {\Lemma (a)} 
Let $\{X_t, t\ge 0\}$ be a random walk on $\Bbb
Z^d$ with transition probability $p(\cdot , \cdot)$ satisfying the
assumptions stated in the beginning of the article. Then, there exists
a finite constant $C$ such that
$$
\aligned
& \sum_{x\in\Bbb Z^d} |p_t(0,x+e_i) - p_t(0,x)|\; \le \; C(1+t)^{-1/2}
\quad \text{and} \\
&\qquad \sum_{x\in\Bbb Z^d} |p_t(0,x+e_i) - p_t(0,x)|^2 \;
\le \; C(1+t)^{-(1+d)/2}\; .
\endaligned
$$
Here $\{e_i, 1\le i\le d\}$ stands for the canonical basis of $\Bbb
R^d$.
\endproclaim

\demo{Proof}
This result can be proved by a coupling argument indicated to us by
E. Andjel or through the local central limit theorem. With slightly
stronger assumptions on the moments of $p(\cdot)$, the local central
limit theorem gives better estimates (of type $t^{-(d+2)/2}$) for the
second term. In dimension $1$, with the assumptions of
indecomposability and finite second moments, we have
$$
\sum_{x\in\Bbb Z} |p_t(0,x+1) - p_t(0,x)|^2 \;
\le \; C(1+t)^{-(1+\delta)}
\Eq(b)
$$
for some $\delta>0$.

The coupling argument is as follows. Consider two random walks
$X^1_t$, $X^2_t$ with transition probability $p(\cdot)$ on $\Bbb Z^d$
and with initial states $0$ and $e_1$. Assume that $p(e_1)>0$.  We
couple these two random walks in the following way. If the process
$(X^1_t, X^2_t)$ is at $(x^1,x^2)$, for $y\neq e_1$ it jumps to
$(x^1+y, x^2+y)$ at rate $p(y)$, it jumps to $(x^1+e_1, x^2)$ at rate
$p(e_1)$ and it jumps to $(x^1, x^2+e_1)$ at rate $p(e_1)$. We proceed
in this way until they meet. From this time on, they jump
together. With this coupling the difference $X^1_t - X^2_t$ is a
nearest-neighbor, symmetric, one-dimensional random walk with
absorption at the origin that starts from $-1$. A well known bound
gives that the probability for this one-dimensional random walk to
have not reached the origin before time $t$ decays as $t^{-1/2}$. This
estimate permits to prove the lemma in the case $p(e_1)>0$. In the
other case, since the transition probability is assumed to be
indecomposable, there exists an integer $n$ and a sequence $0=x_0,
x_1, \dots, x_n=e_1$ such that $p(x_i, x_{i+1})>0$ and we may proceed
in a similar way. Notice that this proof does not require any
assumption on the moments of $p$. \quad \qed
\enddemo

\bigskip
\noindent {\bf Acknowledgments.}

The authors would like to thank the fruitful discussions with
E. Andjel. This research was supported by PRONEX 41.96.0923.00
``Fen\^omenos Cr\'\i ticos em Probabilidade e Processos
Estoc\'asticos'', FAPESP (98/03382-0, AG), FAPERJ (E26/150940/99, CL)
and CNPq grants 301301-79 (AG) and 300358/93-8 (CL)
\bigskip \bigskip

\subheading{References}
\parindent 1cm
\parskip 2mm

\item{[Al]} Aldous D.,
Random walks on finite groups and rapidly mixing Markov chains. Lectures Notes
in Math. {\bf 986} 243-297 (1983). 

\item{[An1]} Andjel E.D.,
A correlation inequality for the symmetric exclusion process. Ann. of
Prob. {\bf 16} 717-721 (1988). 

\item{[An2]} Andjel E.D.,
Finite exclusion process and independent random walk I. Preprint (1994). 

\item{[BG]} Bertein F. and Galves A.,
Comportement asymptotique de deux marches al\'eatoires sur $\Z$
qui int\'eragissent par exclusion. C. R. Acad. Sci. Paris
A {\bf 285} 681-683 (1977).

\item{[BZ]} Bertini, L., Zegarlinski, B., Coercive inequalities for
Kawasaki dynamics:  The product case,  University of Texas Mathematical
Physics archive preprint 96-561. See also Coercive inequalities for
Gibbs measures,  University of Texas Mathematical
Physics archive preprint 96-562.

\item{[CG]} Cancrini N., Galves  A., Approach to equilibrium
in the symmetric simple exclusion process. Markov Proc. Rel. Fields {\bf 2}
175--184 (1995).

\item{DFIP]} De Masi A., Ferrari P.A., Ianiro N., Presutti E.,
Small deviations from local equilibrium for a process which
exhibits hydrodynamical behavior.
J. Stat. Phys.,  {\bf 29} 81-94 (1982).

\item{[DP]} De Masi A. and Presutti E.,
Probability estimates for symmetric simple exclusion random walks.
Ann. Inst. H. Poincar\'e,  Sect. B {\bf 19}, 71-85 (1983).

\item{[De]} Deuschel J.-D.,
Algebraic $L^{2}$ decay of attractive critical processes on the 
lattice. Ann. Probab. {\bf 22}, 264-283 (1994).

\item{[Di]} Diaconis P.,
Group representations in probability and statistics.
IMS Lecture Series {\bf 11} Hayward CA, (1988).

\item{[DS]} Diaconis P., Stroock D.,
Geometric bounds for eigenvalues of Markov chains.
Ann. Appl. Probab. {\bf 1}, 36-61 (1991).

\item{[FPSV]}
Ferrari P.A., Presutti E., Scacciatelli E., Vares M.E.,
The symmetric simple exclusion process, I: Probability estimates.
Stochastic Process. Appl. {\bf 39}, 89-105 (1991).

\item{[H]}
Harris T.E., Additive set valued Markov processes and graphical methods.
Ann. Probab. {\bf 6}, 355-378 (1978).

\item{[GG]} Guiol H., Galves  A.: Relaxation time of the one dimensional
zero range process with constant rate. Markov Proc. Rel. Fields {\bf 3},
323--332, (1997).

\item{[HR]} Hoffman, J.R., Rosenthal, J.S., Convergence of independent
particle systems, Stoch. Proc. Appl. {\bf 56}, 295--305, (1995).

\item{[H]} Holley R.,
Rapid convergence to equilibrium in ferromagnetic stochastic Ising models.
Resenhas IME-USP {\bf 1}, 131-149 (1993).

\item{[JLQY]} Janvresse, E., Landim C., Quastel J. and Yau, H.T.;
Relaxation to equilibrium of conservative dynamics I:
Zero Range processes. Ann. Probab. {\bf 27}, 325--360, (1999).

\item{[L]} Landim C., Decay to equilibrium in $L^\infty$ of finite
interacting particle systems in infinite volume. Markov
Proc. Rel. Fields {\bf 4}, 517--534, (1998).

\item{[Li]} Liggett T.M., {\it Interacting Particle Systems.}
Springer, Berlin (1985).

\item{[P]}
Petrov, V. V., {\it Sums of Independent Random Variables}
Springer Verlag, New York (1975). 

\item{[Q]} Quastel J.,
Diffusion of color in the simple exclusion process.
Comm. Pure and App. Math., {\bf 45} 623-679 (1992).

\item{[S1]} Spitzer F. : Interaction of Markov processes. Adv. Math., {\bf 5}
246-290 (1970). 

\item{[S2]} Spitzer F.,
Recurrent random walk of an infinite particle system.
Trans. Amer. Math. Soc., {\bf 198} 191-199 (1974).

\enddocument